\documentclass{article}
\usepackage{graphicx} 
\usepackage{amsmath}
\usepackage{amsfonts}
\usepackage{amssymb}
\usepackage{dsfont}
\usepackage{amsthm}
\usepackage{xcolor}
\usepackage{authblk}
\usepackage{stmaryrd}

\usepackage{geometry}
\usepackage[
    backend=biber, 
    style=authoryear-ibid, 
    style=numeric, 
    maxcitenames=2, 
    maxbibnames=6,
    sorting=nyt
]{biblatex}

\newtheorem{thm}{Théorèm}[section]
\newtheorem{Corollary}[thm]{Corollary}

\newtheorem{Proposition}[thm]{Proposition}

\newtheorem{Theorem}[thm]{Theorem}
\newtheorem{Definition}[thm]{Definition}

\newtheorem{Remark}[thm]{Remark}

\title{A Constructive Approach for Building Wavelet Bases in \( L^2(\mathbb{R}^d, \mathbb{R}^m) \) with Optimal Properties}

\author{Hicham TARIF, Nadir MAAROUFI}

\date{March 2024}

\begin{document}

\maketitle

\begin{abstract}
The main contribution of this paper is a constructive method for building separable multivariate vector-valued wavelet bases in the general framework of \( L^2(\mathbb{R}^d, \mathbb{R}^m) \) for any \( d, m \geq 1 \). While separable wavelet bases in \( L^2(\mathbb{R}^d, \mathbb{R}) \) are well-established and widely applied, the explicit construction of truly vector-valued wavelet bases remains an open problem, even in the simplest case of \( L^2(\mathbb{R}, \mathbb{R}^2) \), let alone in \( L^2(\mathbb{R}^2, \mathbb{R}^2) \). In practice, the conventional approach applies standard separable wavelet bases of \( L^2(\mathbb{R}^2, \mathbb{R}) \) independently to each component of vector-valued signals in \( L^2(\mathbb{R}^2, \mathbb{R}^2) \). However, this approach fails to capture the intrinsic vectorial structure of the signals. To address this limitation, we propose a constructive approach within the vector-valued wavelet framework, providing a systematic method for constructing such bases in the general case of \( L^2(\mathbb{R}^d, \mathbb{R}^m) \). By linking \( m \)-multiwavelets to vector-valued wavelets, our approach not only enables the systematic construction of separable multivariate bases in \( L^2(\mathbb{R}^d, \mathbb{R}^m) \) that satisfy the vector-valued multiresolution analysis but also ensures that these bases inherit key structural properties, making them well-suited for practical applications.
\end{abstract}


\section{Introduction}

Since the 1990s, wavelets have gained significant attention in signal processing and data analysis due to their ability to capture both frequency and time-domain information simultaneously. This dual capability makes them particularly well-suited for analyzing non-stationary signals. Wavelets have become a standard theoretical framework for analyzing signals \cite{baldazzi2020systematic}, processes \cite{gogolewski2020influence}, and multifractal functions \cite{jaffard2004wavelet}. A fundamental concept in this framework is multiresolution decomposition (MRA), which enables the representation of a function at different levels of resolution. By capturing local and global information, wavelets provide a flexible tool for tasks ranging from function approximation to data compression and reconstruction from samples \cite{alaifari2017reconstructing}. This versatility has revolutionized the analysis of functions and signals in theoretical mathematics and engineering applications. The origins of wavelets date back to 1909 when Alfred Haar introduced the first wavelet basis \cite{lepik2014haar}, long before the term "wavelet" was coined. The Haar wavelet, while simple and orthogonal, is discontinuous, limiting its applicability in some practical contexts. However, Haar wavelets have proven useful in various theoretical contexts \cite{lepik2005numerical}.

Significant advances in wavelet theory were achieved in the 1980s, notably by Yves Meyer \cite{mallat1999wavelet} and Ingrid Daubechies \cite{daubechies1988orthonormal}. Meyer's wavelets are constructed using a band-limited function, ensuring smoothness and regularity in the time domain. Although they lack compact support, their smoothness and orthogonality make them suitable for applications requiring continuously differentiable wavelets. In contrast, Daubechies wavelets are constructed to possess compact support and a specified number of vanishing moments. These properties enhance localization and facilitate the approximation of smooth functions. A key advantage of Daubechies wavelets is their ability to combine compact support with adjustable vanishing moments, making them particularly effective for smooth function approximation, signal analysis, and image processing \cite{daubechies1988orthonormal,zaynidinov2023application, yusuf2020analysis}. The vanishing moments of the Daubechies wavelets play a crucial role in the multifractal formalism. More precisely, they make the wavelets orthogonal to low-order polynomials, enabling the removal of such trends. Sequentially, wavelet transforms based on Daubechies wavelets characterize H\"older spaces and pointwise regularity. They allow estimating the H\"older exponent through decay conditions in wavelet leaders \cite{barral2023frisch1}.

Yves Meyer's pioneering work on wavelets laid the foundation for Stéphane Mallat's multiresolution analysis (MRA) in the scalar case. MRA provides a systematic and theoretical framework for constructing wavelet bases for spaces such as $L^2(\mathbb{R})$. A key starting point in MRA is to find a scaling function (father wavelet) that satisfies specific properties, from which the corresponding wavelet function (mother wavelet) is derived \cite{mallat1999wavelet}. This framework enables the decomposition of signals across multiple scales, facilitating analysis at various levels of resolution, and permitting the construction of a wavelet basis. MRA also provides a theoretical framework for constructing multivariate scalar wavelet bases in spaces such as $L^2(\mathbb{R}^d, \mathbb{R})$. However, despite the general scheme of MRA, building such multivariate wavelet bases in practice remains challenging \cite{krivoshein2016multivariate}. The primary difficulty lies in finding a suitable multivariate scaling function, and even assuming its existence, characterizing the corresponding mother wavelet is complex due to the involvement of multivariate trigonometric polynomials. A common approach to address this issue is using tensor product wavelets (separable wavelets), built from one-dimensional wavelets. 


The success of scalar wavelet bases in analyzing functions has driven the need for tools capable of handling vector-valued functions, which arise in fields such as signal processing, hydrodynamic turbulence, and vector time series analysis. In response, Xiang-Gen Xia and Bruce W. Suter introduced vector-valued multiresolution analysis (VMRA) \cite{xia1996vector}, extending Mallat's MRA to higher dimensions (in the codomain). VMRA provides a decomposition of spaces $L^2(\mathbb{R}, \mathbb{C}^{m \times m})$ and $L^2(\mathbb{R}, \mathbb{C}^{m})$, where $m \in \mathbb{N}$, thus forming the theoretical basis for constructing vector-valued wavelets. Building on this framework, the authors of \cite{xia1996vector} developed vector transforms paving the way for further developments in constructing orthogonal vector-valued wavelets \cite{chen2007study} and biorthogonal vector-valued wavelets \cite{chen2008biorthogonal, chen2008biorthogonality}. Armed with this VMRA theoretical basis, the authors in \cite{xia1996vector} were able to construct vector-valued Meyer wavelets for the space $L^2(\mathbb{R}, \mathbb{R}^m)$. However, the authors could not give an explicit form of this wavelet. To be more precise, unlike the scalar case, the scaling and mother functions are expressed by an infinite product of matrices (see subsection \ref{prob} for a detailed discussion). To overcome this difficulty, the authors were forced to introduce additional assumptions ensuring the existence of a scaling function that satisfies the VMRA conditions. Another work leveraging VMRA to construct a wavelet basis in \( L^2(\mathbb{R}, \mathbb{R}^m) \) with compact support can be found in \cite{chen2007study}. In the same spirit, the authors assumed the existence of a scaling function that satisfies the VMRA.

Despite these theoretical advances, the practical application of vector-valued wavelets remains largely underdeveloped. In most scientific and engineering fields, vector-valued functions are still analyzed component-wise using classical scalar wavelet bases in $L^2(\mathbb{R}^d, \mathbb{R})$. That is, given a vector function $f: \mathbb{R}^d \to \mathbb{R}^m$ for any \( d \geq 1\) and \(m \geq 2 \), the conventional approach consists of applying a separable scalar wavelet basis independently to each component $f_i$, leading to a separate decomposition for each dimension (see Subsection \ref{component}). While this method is computationally convenient and benefits from the well-established theory of scalar wavelets, it fails to capture the intrinsic interactions between the components of $f$, as highlighted in \cite{xia1996vector}. This limitation is particularly problematic in applications such as fluid dynamics, where the velocity components of a flow field are inherently coupled, or in multi-channel signal processing, where correlations between different channels contain crucial information. Given these challenges, the construction of explicit, structured vector-valued wavelet bases that fully leverage the VMRA framework remains an open problem. 

As detailed below, the construction of vector wavelets in \( L^2(\mathbb{R}, \mathbb{R}^m) \) presents challenges that necessitate a departure from the conventional VMRA approach. The traditional method of constructing wavelet bases begins with the VMRA framework, derives a refinement equation, constructs a scaling function, and subsequently determines a mother wavelet. However, obtaining the vector mother wavelet from the vector scaling function within the VMRA framework introduces significant challenges (see Subsection \ref{prob}). To overcome these difficulties, we propose an alternative approach. The core principle of our method is to directly provide both the vector scaling function and the vector mother wavelet from any univariate scalar wavelet basis in a structured manner. This strategy enables the direct construction of a vector-valued wavelet basis, followed by the verification of the VMRA conditions. By leveraging the \(\ast\)-product and the well-established framework of scalar univariate MRA, this article presents a general construction of vector wavelet bases for \( L^2(\mathbb{R}, \mathbb{R}^m) \). Furthermore, this approach ensures that the resulting vector wavelets retain the key properties of the underlying scalar wavelets, including compact support, vanishing moments, and smoothness.

Several studies have extended the concept of VMRA to the multivariate setting, particularly for \( L^2(\mathbb{R}^d, \mathbb{R}^m) \) (see \cite{chen2009properties, chen2008biorthogonal, wang2009characterization}). Despite these efforts, the construction of multivariate vector wavelets based on VMRA remains an open problem, even in the case of \( L^2(\mathbb{R}^2, \mathbb{R}^2) \). The growing demand for multivariate vector-valued wavelets, such as Daubechies and Meyer wavelets, across various application domains highlights the need for further research. In particular, these wavelets play a crucial role in applications such as fluid dynamics in oceanography and aerodynamics, as well as in vector-valued signal processing \cite{cui2008m, fowler2002wavelet}. To address this challenge, we propose a novel approach for constructing multivariate vector-valued wavelets in \( L^2(\mathbb{R}^d, \mathbb{R}^m) \) using univariate vector-valued wavelets in \( L^2(\mathbb{R}, \mathbb{R}^m) \). More precisely, our method employs a tensor product of \( m \)-multiwavelets in \( L^2(\mathbb{R}) \) to construct a multivariate wavelet basis in \( L^2(\mathbb{R}^d, \mathbb{R}) \). Taking advantage of the connection between vector-valued wavelets in \( L^2(\mathbb{R}, \mathbb{R}^m) \) and \( m \)-multiwavelets, we further develop a systematic framework for constructing separable wavelet bases in \( L^2(\mathbb{R}^d, \mathbb{R}^m) \).  

The paper is organized as follows. Section 2 provides a preliminary overview of vector-valued wavelets and the associated challenges, laying the foundation for subsequent developments. Section 3 examines key properties of the \( * \)-product and introduces our construction of vector-valued wavelets in \( L^2(\mathbb{R}, \mathbb{R}^m) \). Section 4 extends Meyer’s tensor product method to \( m \)-multiwavelets, facilitating the construction of separable bases in \( L^2(\mathbb{R}^d, \mathbb{R}) \). Section 5 presents the construction of multivariate vector-valued wavelets in \( L^2(\mathbb{R}^d, \mathbb{R}^m) \) from scalar wavelets. Finally, Section 6 summarizes the main results and provides directions for future research.

\section{Preliminary}

\subsection{Vector-valued multiresolution analysis}

In this section, we recall the basics of the theory of vector-valued multiresolution decomposition developed in \cite{chen2007study, chen2008biorthogonal, chen2008biorthogonality, xia1996vector}. First of all, let us fix the basic notations and definition of an orthonormal sequence:

\noindent $L^2(\mathbb{R}^d,\mathbb{R}^m)=\{h(x)=(h_1(x),\dots ,h_m(x))^T\;: \int_{\mathbb{R}^d}\sum_{i=1}^{m}|h_i(x)|^2dx <+\infty\},$
the inner product in $L^2(\mathbb{R}^d,\mathbb{R}^m)$ is $$\langle f,g \rangle=\int_{\mathbb{R}^d}\sum_{i=1}^{m}f_i(x)g_i(x)dx,\quad 
\text{and the $*$-product is}\;\;\;
\langle f,g\rangle_*=\int_{\mathbb{R}^d}f(x)g(x)^Tdx.$$
Where the superscript \( T \) means the transpose of the vector. We will drive closely to some important properties of the $*$-product in the section \ref{S1}. $I_{m\times m}$ and $0_{m\times m}$ denote the unite matrix and the null matrix of $\mathcal{M}_{m\times m}(\mathbb{R})$ respectively, we take for norm in $\mathcal{M}_{m\times m}(\mathbb{R})$ (or simply $\mathcal{M}_{m\times m}$) the norm $$\|A\|_{\mathcal{M}_{m\times m}}=\|A\|_{1}=\max_j\sum_i|a_{i,j}|, \quad \text{with}\quad A=(a_{i,j})_{1\leq i,j\leq m}.$$

\noindent For \( f \in L^2(\mathbb{R}^d, \mathbb{R}^m) \), its integration and Fourier transform are defined, respectively, as follows:

\[
\int_{\mathbb{R}^d} f(x) \, dx := 
\begin{pmatrix}
\int_{\mathbb{R}^d} f_1(x) \, dx \\
\int_{\mathbb{R}^d} f_2(x) \, dx \\
\vdots \\
\int_{\mathbb{R}^d} f_m(x) \, dx
\end{pmatrix}, \;\;\;
 \widehat{f}(\xi) := \int_{\mathbb{R}^d} f(x) \exp(-i \langle \xi, x \rangle) \, dx = 
\begin{pmatrix}
\int_{\mathbb{R}^d} f_1(x) \exp(-i \langle \xi, x \rangle) \, dx \\
\int_{\mathbb{R}^d} f_2(x) \exp(-i \langle \xi, x \rangle) \, dx \\
\vdots \\
\int_{\mathbb{R}^d} f_m(x) \exp(-i \langle \xi, x \rangle) \, dx
\end{pmatrix},
\]

\noindent where $ \langle \xi, x \rangle $ denotes the inner product of the real vector $ \xi $ and $ x $.

\noindent In this work, our results are presented within the \( L^2(\mathbb{R}^d,\mathbb{R}^m) \) framework. For \( m = 1 \), we simply write \( L^2(\mathbb{R}^d) \) instead of \( L^2(\mathbb{R}^d, \mathbb{R}) \).


\begin{Definition}
 A sequence of vector-valued functions \( \{U_k(x)\}_{k \in \mathbb{Z}^d}  \subset L^2(\mathbb{R}^d, \mathbb{R}^m) \) is called an orthonormal system if it satisfies
\[
\langle U_k, U_u \rangle_* = \delta_{k,u} I_{m\times m}, \quad k, u \in \mathbb{Z}^d,
\]
where \( \delta_{k,u} \) denotes the generalized Kronecker symbol, i.e., for \( k = u \), \( \delta_{k,u} = 1 \); otherwise, \( \delta_{k,u} = 0 \).
We say that \( U(x) \) is an orthonormal vector-valued function if the system of its translations \( \{U(x - k), k \in \mathbb{Z}^d\} \) is orthonormal in \( L^2(\mathbb{R}^d, \mathbb{R}^m) \).
\end{Definition}

\begin{Definition}
An orthonormal system \( \{U_k(x)\}_{k \in \mathbb{Z}^d}  \subset L^2(\mathbb{R}^d, \mathbb{R}^m) \) is called orthonormal basis of a subspace $V\subset L^2(\mathbb{R}^d, \mathbb{R}^m)$, if for every $f\in V$ there exists a unique $ m \times m $ matrix sequence $ \{Q_k\}_{k \in \mathbb{Z}^d}\in \ell^2(\mathbb{Z}^d)^{m \times m} $ such that
$$
f(x) = \sum_{k \in \mathbb{Z}^d} Q_k U_k(x), \quad x \in \mathbb{R}^d,
$$

\noindent where $\ell^2(\mathbb{Z}^d)^{m \times m} = \left\{ B : \mathbb{Z}^d \to \mathbb{R}^{m \times m}, \; \| B \|^2 := \sum_{j,l=1}^d \sum_{k \in \mathbb{Z}^d} |b_{j,l}(k)|^2 < +\infty \right\}$.

\end{Definition}

\noindent The following definition corresponds to the multiresolution decomposition of the space $L^2(\mathbb{R}^d,\mathbb{R}^m)$ using a $d\times d$ dilation matrix $A$ whose all entries are integers and all its eigenvalues are large than one (see \cite{chen2008biorthogonal}).

\begin{Definition}\label{DMRA}
    We say that $\Phi$ is a scaling function for a vector-valued multiresolution analysis $\{ V_j\} _{j\in \mathbb{Z}}$ of $L^2(\mathbb{R}^d,\mathbb{R}^m) $, if the sequence $\{ V_j\} _{j\in \mathbb{Z}}$ satisfies :

\begin{itemize}
    \item $V_{j}\subset V_{j+1}$, for all $j\in \mathbb{Z}$.
    
    \item $\bigcap_{j\in\mathbb{Z}} V_j=\{0\}$, $\bigcup_{j\in \mathbb{Z}} V_j$ is dense in $ L^2(\mathbb{R}^d,\mathbb{R}^m)$.
    
    \item $h(\cdot)\in V_j$ if and only if $h(A\cdot)\in V_{j+1}$, for all $j\in\mathbb{Z}$.

    \item $h(\cdot)\in V_j$ if and only if $h(\cdot-k)\in V_{j}$, for all $j,k\in\mathbb{Z}$.
    
    \item The sequence $\{\Phi (\cdot-k)\;:\;k\in\mathbb{Z}^d \}$ is an orthonormal basis of $V_0$.
\end{itemize}
\end{Definition}

\begin{Remark}
      This definition represents the general form of the multiresolution decomposition of the space \( L^2(\mathbb{R}^d,\mathbb{R}^m) \), as other cases can be expressed in the form of Definition \ref{DMRA} (by other cases, we refer to $N$-scale multiresolution decompositions). For example, the multiresolution decomposition of the space \( L^2(\mathbb{R}^d,\mathbb{R}^m) \) using the scale \( N \) \cite{chen2009design, chen2009existence, chen2009properties} corresponds to a multiresolution decomposition with the dilation matrix \( A = NI_{d\times d} \).

\end{Remark}

\subsection{Challenges in constructing vector-valued wavelets from VMRA}\label{prob}

The most widely used method for building wavelets is multiresolution analysis (MRA) of the space $L^2(\mathbb{R})$, in combination with Fourier theory. MRA provides a systematic and robust framework for constructing standard wavelet bases. This methodology enables the approximation of $L^2(\mathbb{R})$ at multiple levels of resolution, facilitating the development of wavelet bases with desirable properties, such as those exhibited by Meyer and Daubechies wavelets. This approach has been generalized to the vector-valued function space $L^2(\mathbb{R}, \mathbb{R}^m)$ \cite{chui1996study}. In this section, we review key techniques and properties for constructing vector-valued wavelets in $L^2(\mathbb{R}, \mathbb{R}^m)$. We also examine the challenges inherent in employing VMRA to construct a scaling function and replicate the properties of Daubechies wavelets.


\noindent Consider a scaling function $\Phi$ for the vector-valued multiresolution analysis $\{ V_j\} _{j\in \mathbb{N}}$ of $L^2(\mathbb{R},\mathbb{R}^m) $ (for simplicity we take $A=2I_{1\times 1}=2$). Since \( \Phi \in V_0 \subset V_1 \), there exists a sequence of \( m \times m \) matrix \( \{P_k\}_{k \in \mathbb{Z}} \in \ell^2(\mathbb{Z})^{m \times m} \), which has finite non-zero terms (because we want a compactly supported scaling function), such that
\begin{equation}\label{E1}
    \Phi(x) = \sum_{k \in \mathbb{Z}} P_k \Phi(2x - k).
\end{equation}

\noindent Equation \ref{E1} is called a \textit{refinement} equation. Taking the Fourier transform of (\ref{E1}), we get
\begin{equation}
    \widehat{\Phi}(\xi) = P\left(\frac{\xi}{2}\right) \widehat{\Phi}\left(\frac{\xi}{2}\right), \quad \xi \in \mathbb{R},
\end{equation}

\noindent where
\begin{equation}\label{E2}
    P(\xi) = \frac{1}{2} \sum_{k \in \mathbb{Z}} P_k e^{-i \langle k, \xi \rangle}, \quad \xi \in \mathbb{R},
\end{equation}


\vspace{0.5cm}
\noindent then, as usual in this framework, by assuming that $\widehat{\Phi}$ is continuous at $0$, we get 

$$
\widehat{\Phi}(\xi) = \lim_{k \to \infty}  
P\left(\frac{\xi}{2}\right) 
P\left(\frac{\xi}{4}\right) 
\cdots P\left(\frac{\xi}{2^k}\right)\widehat{\Phi}\left(\frac{\xi}{2^k}\right)
= 
\prod_{k=1}^{+\infty} 
P\left(\frac{\xi}{2^k}\right) \widehat{\Phi}(0).
$$

\noindent Without loss of generality, for simplicity, we can suppose that $\widehat{\Phi}(0)=I_{m\times m}$.  
Using, the orthonormality of the scaling function we straightforwardly obtain the condition: $$ P(\xi) P(\xi)^* + P(\xi + \pi) P(\xi + \pi)^* = I_{m\times m}, \quad \xi \in \mathbb{R},$$
where $*$ means the complex conjugate and the transpose.
These equations are the starting point for constructing a scaling function that satisfies the VMRA. However, due to their matrix form, these equations are required to fulfill some conjectured conditions. More precisely, at this level, authors assumed a formal estimation of the infinite product of the sequence \( \left(P\left( \frac{\xi}{2^n}\right)\right)_{n \geq 1} \) to establish the existence of such a scaling function (see \cite{xia1996vector}). Once a suitable scaling function is constructed, a wavelet mother can be built based on other equations mentioned below (see \cite{chen2007study}).

\vspace{0.5cm}

\noindent Let \( W_j, j \in \mathbb{Z} \) be the orthonormal complementary sub-spaces of \( V_j \) in \( V_{j+1} \), i.e.,
\[
W_j = V_{j+1} \ominus V_j, \quad j \in \mathbb{Z}.
\]
We aim to construct a mother vector-valued wavelet $\Psi$ as an orthonormal basis of $W_0$. Since \( \Psi \in W_0 \subset V_1 \), there exists a unique finitely supported sequence \( \{B_k\}_{k \in \mathbb{Z}} \) of \( m \times m \) matrices such that
\begin{equation}\label{E4}
    \Psi(x) = \sum_{k \in \mathbb{Z}} B_k \Phi(2x - k).
\end{equation}

\noindent Let
\begin{equation}\label{E3}
B(\xi) = \frac{1}{2} \sum_{k \in \mathbb{Z}} B_k \exp(-ik\xi).
\end{equation}
Then, the equation (\ref{E4}) becomes
\begin{equation}
\widehat{\Psi}(\xi) = B\left(\frac{\xi}{2}\right) \widehat{\Phi}\left(\frac{\xi}{2}\right), \quad \xi \in \mathbb{R}.
\end{equation}

\noindent An interesting result in \cite{chen2007study} is the following theorem, which establishes a necessary and sufficient condition for characterizing the mother vector-valued wavelet, assuming the existence of a vector-valued scaling function.

\begin{Theorem}\cite{chen2007study}
Let \( \Phi(\cdot) \) be an orthogonal vector-valued scaling function. Assume \( \Phi(\cdot) \in L^2(\mathbb{R}, \mathbb{R}^m) \), and \( P(\xi) \) and \( B(\xi) \) are defined by \ref{E2} and \ref{E3}, respectively. Then \( \Psi(\cdot) \) is an orthogonal vector-valued wavelet function associated with \( \Phi(\cdot) \) if and only if
\[
P(\xi) B(\xi)^* + P(\xi + \pi) B(\xi + \pi)^* = 0_{m\times m}, \quad \xi \in \mathbb{R},
\]
\[
B(\xi) B(\xi)^* + B(\xi + \pi) B(\xi + \pi)^* = I_{m\times m}, \quad \xi \in \mathbb{R},
\]
where $*$ means the complex conjugate and the transpose.
\end{Theorem}

\noindent Even with this theorem in hand, constructing a concrete vector-valued wavelet basis remains a challenging task. More specifically, this difficulty arises from the lack of key results on polynomial matrices, such as those found in Corollary 5.5.4 and Lemma 6.1.3 in \cite{daubechies1992ten} for the scalar case \( L^2(\mathbb{R}) \). These propositions, central to the construction of wavelets like Daubechies wavelets in \( L^2(\mathbb{R}) \), concern univariate polynomials and do not directly extend to the vector-valued setting. The main obstacle lies in the absence of desirable properties for trigonometric matrices, which prevents us from following the scalar construction process. In Section~\ref{S3}, we propose a novel method for constructing vector-valued wavelets in \( L^2(\mathbb{R},\mathbb{R}^m) \), \( m \geq 1 \), which does not rely on the VMRA framework. Our method retains crucial properties such as compact support, regularity, and a prescribed number of vanishing moments, building upon wavelets in the scalar case.

\subsection{Multiwavelets vs. vector-valued wavelets}\label{component}

To avoid confusion over terminology, let us start by clarifying certain terms used in this article.
\vspace*{-0.4cm}
\paragraph{\( m \)-Multiwavelet:} Throughout this article, the term \textit{\( m \)-multiwavelet} refers to a basis of \( L^2(\mathbb{R}) \) consisting of \( m \) scaling functions \( \phi_1, \dots, \phi_m \) and \( m \) mother wavelets \( \psi_1, \dots, \psi_m \). These functions form an orthonormal system in \( L^2(\mathbb{R}) \) and generate a basis for \( L^2(\mathbb{R}) \) through translations of the \( m \) scaling functions, as well as dilations and translations of the \( m \) mother wavelets. The first example of a 2-multiwavelet in \( L^2(\mathbb{R}) \) was introduced in \cite{geronimo1994fractal}, where the authors constructed two scaling functions \((\phi_1, \phi_2)\) and two mother wavelets \((\psi_1, \psi_2)\) that generate \( L^2(\mathbb{R}) \). This terminology naturally extends to the vector-valued setting \( L^2(\mathbb{R}^d,\mathbb{R}^m) \), where such constructions are referred to as vector-valued multiwavelets in \cite{cui2008m}.
\vspace*{-0.4cm}
\paragraph{Multiwavelets:}  
In this article, the term \textit{multiwavelets} in a space \( L^2(\mathbb{R}^d,\mathbb{R}^m) \) refers to a component-wise approach to vector-valued wavelets. More precisely, multiwavelets in \( L^2(\mathbb{R}^d,\mathbb{R}^m) \) are multivariate vector-valued functions whose components are multivariate scalar scaling functions and mother wavelets in \( L^2(\mathbb{R}^d,\mathbb{R}) \). However, they do not necessarily form an orthonormal basis of \( L^2(\mathbb{R}^d,\mathbb{R}^m) \). For example, multiwavelets in \( L^2(\mathbb{R}^2, \mathbb{R}^2) \) can be constructed from standard separable wavelet bases of \( L^2(\mathbb{R}^2, \mathbb{R}) \), allowing the signal to be analyzed component by component under the implicit assumption of independence between components. Nevertheless, they provide a flexible framework for capturing signal characteristics in various applications (see, for instance, \cite{kestener2004generalizing, sakran2023new}).

In fact, the term multiwavelets first appeared in the scalar case in \cite{geronimo1994fractal}, where the authors constructed two mother wavelets with short support, continuity, and specific symmetry properties, forming an orthonormal basis for \( L^2(\mathbb{R}) \). This innovation enabled the simultaneous use of multiple wavelet functions to analyze scalar signals, offering improved approximation capabilities for signals with complex structures. This concept was later extended to \( L^2(\mathbb{R}, \mathbb{R}^m) \) with a different perspective, where multiwavelets are used to analyze the components of vector-valued signals simultaneously using a vector of scalar wavelet functions. In this case, the analysis consists of selecting a scalar wavelet basis and processing each component of the vector separately. However, a fundamental limitation of the multiwavelet approach is that it decomposes each vector component independently, capturing only intra-component correlations in the time or spatial domain while failing to explicitly account for inter-component dependencies \cite{cui2008m, fowler2002wavelet}.

The gap in the multiwavelet framework was first highlighted by Xia et al. in \cite{xia1996vector}, who introduced a more comprehensive theoretical framework extending the wavelet transform to the vector case. Building on this, their work developed the concept of vector-valued wavelets, aimed at decomposing vector-valued signals while taking into account relationships between vector components. The vector-valued wavelet approach relies on the VMRA, enabling simultaneous treatment of correlations in time and across vector components. As noted in \cite{xia1996vector}, vector-valued wavelets are designed to decorrelate vector-valued signals not only in the time domain but also across vector components at a fixed time point. This inter-component decorrelation makes vector-valued wavelets particularly effective for applications where dependencies between components are critical, such as multivariate time series and color image processing. In contrast, multiwavelets primarily focus on time-domain decorrelation without an intrinsic mechanism to address inter-component correlations, making them less suitable for applications where inter-component relationships play a significant role. In summary, vector-valued wavelets constitute a broader framework than multiwavelets.

One of the main contributions of \cite{xia1996vector} is the establishment of a connection between vector-valued wavelets and \( m \)-multiwavelets, presented within the general context of \( L^2(\mathbb{R},\mathbb{C}^m) \). This framework is also directly applicable to the case of \( L^2(\mathbb{R},\mathbb{R}^m) \).

\begin{Proposition}\label{XP2}
Let \( \Phi = (\phi_1, \dots, \phi_m) \) be a vector-valued scaling function associated with a VMRA of \( L^2(\mathbb{R},\mathbb{R}^m) \), and let \( \Psi = (\psi_1, \dots, \psi_m) \) be its corresponding vector-valued mother wavelet. Then, the components \( \phi_l \) for \( l = 1, 2, \ldots, m \), constitute \( m \)-scaling functions, and the components \( \psi_l \) for \( l = 1, 2, \ldots, m \), constitute \( m \)-mother wavelets of \(L^2(\mathbb{R})\).
\end{Proposition}

The if and only if version of this proposition is established in Corollary \ref{cor}.

\section{Construction of vector-valued wavelets in $L^2(\mathbb{R},\mathbb{R}^m)$ from scalar wavelets}\label{S3}

The main consequence of Proposition \ref{XP2} is that vector-valued wavelets of \( L^2(\mathbb{R},\mathbb{R}^m) \) represent a broader framework than \( m \)-multiwavelets of \( L^2(\mathbb{R}) \), as they can be generated from their component functions. In this section, we begin by investigating some properties of the \( * \)-product needed in the sequel of this article. Then, we present a method for constructing an orthonormal vector-valued wavelet basis for \( L^2(\mathbb{R},\mathbb{R}^m) \) that respects the VMRA, starting from a given wavelet basis of \( L^2(\mathbb{R}) \). Additionally, we establish the connection between \(m\)-multiwavelet and vector-valued wavelet in \( L^2(\mathbb{R},\mathbb{R}^m) \).

\subsection{The properties of $*$-product in $L^2(\mathbb{R}^d,\mathbb{R}^m)$}\label{S1}
The results are for $d=1$ and $m=2$, but it can be easily extended to the general case.

\begin{Proposition}[Continuity]\label{P1}
    The linear mapping 
\begin{align*}
    \langle f,\cdot\rangle_*:(L^2(\mathbb{R},\mathbb{R}^2),\|\cdot\|_{L^2(\mathbb{R},\mathbb{R}^2)})        &\rightarrow (\mathcal{M}_{2\times 2},\|\cdot\|_{\mathcal{M}_{2\times 2}})\\
    g \quad\quad \quad\quad \quad\quad&\rightarrow  \langle f,g\rangle _*
\end{align*}
is continuous, where $\|(a_{i,j})_{i,j}\|_{\mathcal{M}_{2\times 2}}=\max_j\sum_i|a_{i,j}|$, and $f\in L^2(\mathbb{R},\mathbb{R}^2)$.
\end{Proposition}

\noindent \textbf{Proof:} The linearity of the application is straightforward, let $g=(g_1,g_2)^T\in L^2(\mathbb{R},\mathbb{R}^2)$ we have :

\begin{equation*}
\langle f,g\rangle_*=\int_{\mathbb{R}}\begin{pmatrix}
              f_1(x)g_1(x) &  f_1(x)g_2(x) \\
              f_2(x)g_1(x) & f_2(x)g_2(x)
\end{pmatrix}    dx  ,
\end{equation*}
so 
\begin{align*}
\|\langle f,g\rangle_*\|_{\mathcal{M}_{2\times 2}} &=\max\{ |\int_{\mathbb{R}}f_1(x)g_1(x)dx |+|\int_{\mathbb{R}}f_2(x)g_1(x)dx |, |\int_{\mathbb{R}}f_1(x)g_2(x)dx |+|\int_{\mathbb{R}}f_2(x)g_2(x)dx |     \} \\
            &\leq |\int_{\mathbb{R}}f_1(x)g_1(x)dx |+|\int_{\mathbb{R}}f_2(x)g_1(x)dx |+|\int_{\mathbb{R}}f_1(x)g_2(x)dx |+|\int_{\mathbb{R}}f_2(x)g_2(x)dx |  \\
            &\leq \|f_1\|_{L^2(\mathbb{R})}\|g_1\|_{L^2(\mathbb{R})}+\|f_2\|_{L^2(\mathbb{R})}\|g_1\|_{L^2(\mathbb{R})}+\|f_1\|_{L^2(\mathbb{R})}\|g_2\|_{L^2(\mathbb{R})}+\|f_2\|_{L^2(\mathbb{R})}\|g_2\|_{L^2(\mathbb{R})} \\
            &\leq 2\|f\|_{L^2(\mathbb{R},\mathbb{R}^2)}\|g_1\|_{L^2(\mathbb{R})}+2\|f\|_{L^2(\mathbb{R},\mathbb{R}^2)}\|g_2\|_{L^2(\mathbb{R})} \\
            &\leq 4\|f\|_{L^2(\mathbb{R},\mathbb{R}^2)}\|g\|_{L^2(\mathbb{R},\mathbb{R}^2)}
\end{align*}
because $\|f\|_{L^2(\mathbb{R},\mathbb{R}^2)}^2=\|f_1\|_{L^2(\mathbb{R})}^2+\|f_2\|_{L^2(\mathbb{R})}^2$, so $\|f_1\|_{L^2(\mathbb{R})}+\|f_2\|_{L^2(\mathbb{R})}\leq2\|f\|_{L^2(\mathbb{R},\mathbb{R}^2)}$.

\hfill \qed

\begin{Proposition}[Hadamard product]\label{lem}
Let $h(x,y)=(f_1(x)g_1(y),f_2(x)g_2(y))^T$ , $k(x,y)=(r_1(x)u_1(y),r_2(x)u_2(y))^T\in L^2(\mathbb{R}^2,\mathbb{R}^2)$, where $f,g,r,u \in L^2(\mathbb{R})$.
Then $$\langle h,k\rangle_*=\langle f,r\rangle _*\otimes \langle g,u\rangle_*   ,$$
in which $(a_{i,j})_{i,j\in \{1,\dots ,n\}}\otimes (b_{i,j})_{i,j\in\{1,\dots ,n\}}=(a_{i,j}b_{i,j})_{i,j\in\{1,\dots ,n\}}$.
    
\end{Proposition} 
\noindent \textbf{Proof:} Indeed 
\begin{align*}
    \langle h,k\rangle_* 
    &= \int_{\mathbb{R}^2} h(x,y)k(x,y)^Tdxdy \\[8pt]
    &= \int_{\mathbb{R}^2} 
    \begin{pmatrix}
        f_1(x)g_1(y)r_1(x)u_1(y) & f_1(x)g_1(y)r_2(x)u_2(y) \\[5pt]
        f_2(x)g_2(y)r_1(x)u_1(y) & f_2(x)g_2(y)r_2(x)u_2(y)
    \end{pmatrix} 
    dxdy \\[10pt]
    &=  
    \begin{pmatrix}
        \displaystyle\int_{\mathbb{R}} f_1(x)r_1(x)dx \int_{\mathbb{R}} g_1(y)u_1(y)dy 
        & \displaystyle\int_{\mathbb{R}} f_1(x)r_2(x)dx \int_{\mathbb{R}} g_1(y)u_2(y)dy \\[10pt]
        \displaystyle\int_{\mathbb{R}} f_2(x)r_1(x)dx \int_{\mathbb{R}} g_2(y)u_1(y)dy 
        & \displaystyle\int_{\mathbb{R}} f_2(x)r_2(x)dx \int_{\mathbb{R}} g_2(y)u_2(y)dy
    \end{pmatrix} \\[12pt]
    &=  
    \left( \int_{\mathbb{R}} 
    \begin{pmatrix}
        f_1(x)r_1(x) & f_1(x)r_2(x) \\[5pt]
        f_2(x)r_1(x) & f_2(x)r_2(x)
    \end{pmatrix} 
    dx \right) 
    \otimes  
    \left( \int_{\mathbb{R}} 
    \begin{pmatrix}
        g_1(y)u_1(y) & g_1(y)u_2(y) \\[5pt]
        g_2(y)u_1(y) & g_2(y)u_2(y)
    \end{pmatrix} 
    dy \right) \\[12pt]
    &= \langle f,r\rangle_* \otimes \langle g,u\rangle_*.
\end{align*}

\hfill \qed

\begin{Remark}\label{R1}
    This application $\langle \cdot,\cdot\rangle_*$ has linearity corresponding to the first and second variables, and if $A$ is a matrix we have $$\langle f, Ag\rangle_*=\langle f,g\rangle_*A^T,\quad\text{and}\quad \langle f,g\rangle_*=\langle g,f\rangle_*^T, $$
    where $A^{T}$ is the transpose matrix of $A$. 
\end{Remark}

\subsection{The construction of orthonormal vector-valued wavelets in $L^2(\mathbb{R},\mathbb{R}^m)$}\label{S}

The following proposition demonstrates how to construct a wavelet basis for the space \( L^2(\mathbb{R},\mathbb{R}^m) \) starting from a scalar wavelet basis of the space \( L^2(\mathbb{R}) \).

\begin{Proposition}\label{rm}
Let $\phi$ be a scaling function and $\psi$ a mother wavelet of the space $L^2(\mathbb{R})$ generated from  MRA, then the system of functions 
$$\Phi_{k}(x)=(\phi(x-k),\psi(x-k),\dots ,2^{(m-2)/2}\psi(2^{(m-2)}x-k))^T,$$  and $$\Psi_{j,k}(x)=(2^{(j+m-1)/2}\psi(2^{j+m-1}x-k), 2^{(j+m)/2}\psi(2^{j+m}x-k),\dots,2^{(j+2m-2)/2}\psi(2^{j+2m-2}x-k) )^T, $$
form a wavelet basis of the space \( L^2(\mathbb{R},\mathbb{R}^m) \) with dilation \( 2^m \) that respects the VMRA and inheriting the properties of the scalar wavelet.
    
\end{Proposition}

\noindent\textbf{Proof:} We show the results of this section for the case $m=2$, the proof is similar for the cases where $m>2$. In our case, the dilation matrix $A$ is simply the scalar $4$.

Consider an orthogonal wavelets wavelet $\{\phi_k, \psi_{j,k}\}_{k\in\mathbb{Z},j\geq 0}$ in $L^2(\mathbb{R})$, we have $L^2(\mathbb{R})=V_0\bigoplus_{j\geq 0 } W_j$, where $V_0=\overline{Vect\{\phi(.-k)\;:\;k\in\mathbb{Z}  \}} $ and $W_j=\overline{Vect\{ 2^{j/2}\psi(2^j.-k)\;:\;k\in\mathbb{Z}  \}} $. Let us consider $\{\Phi_k\}_{k\in\mathbb{Z}} $ and $\{\Psi_{j,k}\}_{j\in 2\mathbb{N}+1,k\in\mathbb{Z}}$ such that : $$\Phi_{k}(x)=(\phi(x-k),\psi(x-k))^T \quad \text{and} \quad \Psi_{j,k}(x)=(2^{(j+1)/2}\psi(2^{j+1}x-k), 2^{(j+2)/2}\psi(2^{j+2}x-k) )^T $$
It is essay to see that $\Phi_{k} $ and $\Psi_{j,k} $ are in $L^2(\mathbb{R},\mathbb{R}^2)$. Furthermore, this system $\{\Phi_{k},\Psi_{j,k} \}$ is a wavelet basis of $ L^2(\mathbb{R},\mathbb{R}^2)$ with respect to the scalar dilation $m=4$. 

The system is orthonormal because we have 

\begin{equation*}
\langle\Phi_k,\Phi_{k'}\rangle_*=\int_{\mathbb{R}}\begin{pmatrix}
              \phi(x-k)\phi(x-k') &  \phi(x-k)\psi(x-k') \\
              \psi(x-k)\psi(x-k') & \psi(x-k)\psi(x-k')
\end{pmatrix}    dx  =\delta_{k,k'}I_{2\times 2}, \quad \text{for all $k,k'\in\mathbb{Z}$}
\end{equation*}

\begin{equation*}
 \langle\Phi_k,\Psi_{j,k'}\rangle_*=\int_{\mathbb{R}}\begin{pmatrix}
              2^{(j+1)/2}\phi(x-k)\psi(2^{j+1}x-k') &  2^{(j+2)/2}\phi(x-k)\psi(2^{j+2}x-k') \\
              2^{(j+1)/2}\psi(x-k)\psi(2^{j+1}x-k') & 2^{(j+2)/2}\psi(x-k)\psi(2^{j+2}x-k')
\end{pmatrix}    dx  =0_{2\times 2},
\end{equation*}

because $$\int_{\mathbb{R}}\phi(x-k)\psi(2^{j+1}x-k')dx=0 ,\quad \text{for all $j\geq 0$ and $k,k'\in\mathbb{Z}$}, $$
and $$\int_{\mathbb{R}}\psi(x-k)\psi(2^{j+1}x-k')dx=0 ,\quad \text{for all $j\geq 0$ and $k,k'\in\mathbb{Z}$}. $$

Also  
$$
\langle\Psi_{j,k},\Psi_{j',k'}\rangle_*=\int_{\mathbb{R}}\begin{pmatrix}
              2^{(j+1)/2}\psi(2^{j+1}x-k)2^{(j'+1)/2}\psi(2^{j'+1}x-k') &  2^{(j+1)/2}\psi(2^{j+1}x-k)2^{(j'+2)/2}\psi(2^{j'+2}x-k') \\
              2^{(j+2)/2}\psi(2^{j+2}x-k)2^{j'+2}\psi(2^{j'+1}x-k') & 2^{(j+2)/2}\psi(2^{j+2}x-k)2^{j'+2}\psi(2^{j'+2}x-k')
\end{pmatrix}dx $$
 $$=\delta_{j,j'}\delta_{k,k'}I_{2\times 2}, \quad \text{for all $j,j'\in2\mathbb{N}+1$, and $k,k'\in\mathbb{Z}$}.
$$
So the system $\{\Phi_k\}_{k\in\mathbb{Z}} $ and $\{\Psi_{j,k}\}_{j\in 2\mathbb{N}+1,k\in\mathbb{Z}}$ is orthonormal in $L^2(\mathbb{R},\mathbb{R}^2)$.

\vspace{0.5cm}

The next step is to show that if $f\in L^2(\mathbb{R},\mathbb{R}^2)$ there exist a system of $2\times 2$ size matrices $\{C_k, C_{j,k} \}$, such that $$f=\sum_{k\in\mathbb{Z}}C_k\Phi_k+\sum_{k\in\mathbb{Z},
    j\in 2\mathbb{N}+1
 }C_{j,k}\Psi_{j,k}.$$

Let $f\in L^2(\mathbb{R},\mathbb{R}^2)$, so there is two functions in $L^2(\mathbb{R})$ named $f_1$ and $f_2$ such that $f=(f_1,f_2)^T$

$$f_1=\sum_{k\in\mathbb{Z}}c^1_k\phi_k+\sum_{k\in\mathbb{Z},j\in\mathbb{N}
    }c^1_{j,k}\psi_{j,k},$$

$$f_2=\sum_{k\in\mathbb{Z}}c^2_k\phi_k+\sum_{k\in\mathbb{Z},j\in\mathbb{N}}c^2_{j,k}\psi_{j,k}.$$

\noindent Now, since wavelet bases are unconditional in \(L^2(\mathbb{R})\), thus the reordering of indices or terms does not change the sum. We can rearrange the sum in the following manner.
$$f_1=\sum_{k\in\mathbb{Z}}(c^1_k\phi_k+c^1_{0,k}\psi_{0,k})+\sum_{j\in2\mathbb{N}+1}\sum_{k\in\mathbb{Z}}(c^1_{j,k}\psi_{j,k}+c^1_{j+1,k}\psi_{j+1,k}),$$
$$f_2=\sum_{k\in\mathbb{Z}}(c^2_k\phi_k+c^2_{0,k}\psi_{0,k})+\sum_{j\in2\mathbb{N}+1}\sum_{k\in\mathbb{Z}}(c^2_{j,k}\psi_{j,k}+c^2_{j+1,k}\psi_{j+1,k}).$$

\noindent This implies that there exist a system of $2\times 2$ size matrices $\{C_k, C_{j,k} \}$, such that $$f=\sum_{k\in\mathbb{Z}}C_k\Phi_k+\sum_{k\in\mathbb{Z},
    j\in 2\mathbb{N}+1
 }C_{j,k}\Psi_{j,k}.$$
Where the convergence is in the Hilbert space $(L^2(\mathbb{R},\mathbb{R}^2),\|.\|_{L^2(\mathbb{R},\mathbb{R}^2)}) $. Furthermore, these matrices are :

$$C_k=\begin{pmatrix}
              c^1_k & c^1_{0,k} \\
              c^2_k & c^2_{0,k}
\end{pmatrix} ,$$

$$C_{j,k}=\begin{pmatrix}
              c^1_{j,k} & c^1_{j+1,k} \\
              c^2_{j,k} & c^2_{j+1,k}
\end{pmatrix} .$$

\noindent Hence, by proposition \ref{P1}, remark \ref{R1} and the orthonormality of the system $\{\Psi_k,\Phi_{j,k}\}$, the continuity of the application $\langle.,.\rangle_*$ implies :
\begin{align*}
\langle\sum_{k'\in\mathbb{Z}}C_{k'}\Phi_{k'}+\sum_{k'\in\mathbb{Z},
    j\in 2\mathbb{N}+1
 }C_{j,k'}\Psi_{j,k'},\Psi_{j,k}\rangle_* &=\sum_{k'\in\mathbb{Z}}\langle C_{k'}\Phi_{k'},\Psi_{j,k}\rangle_*+\sum_{k'\in\mathbb{Z},j\in2\mathbb{N}+1}\langle C_{j,k'}\Psi_{j,k'},\Psi_{j,k} \rangle_*\\
  &=\langle C_{j,k}\Psi_{j,k},\Psi_{j,k}\rangle_*\\
  &=C_{j,k}\langle\Psi_{j,k},\Psi_{j,k}\rangle_*^T\\
  &=C_{j,k}.
\end{align*}

\noindent The case of $C_k$ is similar.


\noindent Let us check that our scaling function respects the VMRA. Let $V_0=\overline{Vect\{C\Phi_k \;:\; k\in\mathbb{Z},\;C\in\mathcal{M}_{2\times 2}  \}}$, and  $W_j=\overline{Vect\{C\Psi_{j,k}\;:\;k\in\mathbb{Z},\;C\in\mathcal{M}_{2\times 2} \}}$ where the closer is taken with respect to $L^2(\mathbb{R},\mathbb{R}^2)$ norm.
$V_0$ and $\{V_j\}_{j\in \mathbb{Z}}$ are closed subspaces of $L^2(\mathbb{R},\mathbb{R}^2)$, and since the vector components of $\Phi$ are the scaling function and mother wavelet of the space $ L^2(\mathbb{R})$ we can obtain:

\begin{itemize}
    \item  $V_{j}   \subset V_{j+1}$, for all $j\in \mathbb{Z}$, where $V_{j}=\overline{Vect\{C\Phi(A^j\cdot-k) \;:\;k\in\mathbb{Z},\;C\in\mathcal{M}_{2\times 2} \}}$.
    
    \item $\bigcap_{j\in\mathbb{Z}} V_j=\{0\}$, $\bigcup_{j\in \mathbb{Z}} V_j$ is dense in $ L^2(\mathbb{R},\mathbb{R}^2)$. 
    
    \item $h(\cdot)\in V_j$ if and only if $h(A\cdot)\in V_{j+1}$, for all $j\in\mathbb{Z}$.

    \item $h(\cdot)\in V_j$ if and only if $h(\cdot-k)\in V_{j}$, for all $j,k\in\mathbb{Z}$.
    
    \item The sequence $\{\Phi (\cdot-k)\;:\;k\in\mathbb{Z} \}$ is an orthonormal basis of $V_0$.
\end{itemize}
 Additionally, the space $W_0$ is orthogonal to $V_0$. Therefore we deduce :

$$L^2(\mathbb{R},\mathbb{R}^2)=V_0 \bigoplus_{j\in \mathbb{N} }W_{A^j}, $$

\noindent where $W_{A^j}=\overline{Vect\{C\Psi_k(A^jx)\;:\; k\in\mathbb{Z} ,\;C\in\mathcal{M}_{2\times 2}\}}^{L^2(\mathbb{R},\mathbb{R}^2)}$. Recall that $$\Psi(x)=(\psi(2x),\psi(4x))^T,\;\text{and}\;\Psi_k(A^jx)=\Psi(A^jx-k).$$
\hfill \qed

 \begin{Corollary}\label{cor}
Consider the following wavelet bases:
\begin{enumerate}
    \item A wavelet basis in \(L^2(\mathbb{R})\)  respecting the MRA.
    \item A vector-valued wavelet basis respecting the VMRA in \(L^2(\mathbb{R}, \mathbb{R}^m)\) for \(m \ge 2\).
    \item An \(m\)-multiwavelet basis in \(L^2(\mathbb{R})\) respecting the MRA.
\end{enumerate}
Then
\[
  \textup{(1)} \;\implies\; \textup{(2)}
  \quad\text{and}\quad
  \textup{(2)} \;\iff\; \textup{(3)}.
\]
In particular, from any wavelet basis \textup{(1)} in \(L^2(\mathbb{R})\), one can construct both \textup{(2)} and \textup{(3)}. Moreover, the construction preserves key properties such as regularity, vanishing moments, and compact support.
\end{Corollary}

\noindent\textbf{Proof:} First, (1) $\Rightarrow$ (2) is due to the Proposition \ref{rm}. Second, (2) $\Rightarrow$ (3) is due to Proposition \ref{XP2} from \cite{xia1996vector}. Furthermore, the \(m\)-multiwavelet obtained respect the MRA in the sense that $V_j=\overline{\operatorname{Vect} \{ \phi_i(2^{j}\cdot - k) : i \in \{1, \dots, m\}, k \in \mathbb{Z} \} }$ form a sequence of nested subspaces in \(L^2(\mathbb{R})\).

\noindent For (3) $\Rightarrow$ (2), let \( \phi_1, \dots, \phi_m \)  \( m \)-scaling functions and \( \psi_1, \dots, \psi_m \) \( m \)-mother wavelets, we can construct a vector-valued wavelet basis in the space \( L^2(\mathbb{R}, \mathbb{R}^m) \). In fact, any scalar function \( g \in L^2(\mathbb{R}) \) can be rewritten as
\[
g(x) = \sum_{k \in \mathbb{Z}} \sum_{l=1}^m c_k^l \phi_l(x-k) + \sum_{j \in \mathbb{N}, k \in \mathbb{Z}} \sum_{l=1}^m c_{j,k}^l \psi_l(2^j x-k),
\]
where the coefficients \( \{c_k^l\} \) and \( \{c_{j,k}^l\} \) depend on the function \( g \).

\noindent Now, if we take a vector-valued function \( f := (f_1, \dots, f_m)^T \in L^2(\mathbb{R}, \mathbb{R}^m) \), we have:
\[
\begin{cases}
    f_1(x) = \sum_{k \in \mathbb{Z}} \sum_{l=1}^m c_k^{1,l} \phi_l(x-k) + \sum_{j \in \mathbb{N}, k \in \mathbb{Z}} \sum_{l=1}^m c_{j,k}^{1,l} \psi_l(2^j x-k) \\
    \vdots \\
    f_m(x) = \sum_{k \in \mathbb{Z}} \sum_{l=1}^m c_k^{m,l} \phi_l(x-k) + \sum_{j \in \mathbb{N}, k \in \mathbb{Z}} \sum_{l=1}^m c_{j,k}^{m,l} \psi_l(2^j x-k).
\end{cases}
\]
This leads to the vectorized representation:
\[
f(x) = \sum_{k \in \mathbb{Z}} C_k \Phi(x-k) + \sum_{j \in \mathbb{N}, k \in \mathbb{Z}} C_{j,k} \Psi(2^j x-k),
\]
where $\Phi=(\phi_1,\dots,\phi_m)^T$ and $\Psi=(\psi_1,\dots,\psi_m)^T$. Also, \( C_k \) and \( C_{j,k} \) are \( m \times m \) matrices. Specifically:
\[
C_k =
\begin{pmatrix}
c_k^{1,1} & \cdots & c_k^{1,m} \\
\vdots & \ddots & \vdots \\
c_k^{m,1} & \cdots & c_k^{m,m}
\end{pmatrix}, \quad
C_{j,k} =
\begin{pmatrix}
c_{j,k}^{1,1} & \cdots & c_{j,k}^{1,m} \\
\vdots & \ddots & \vdots \\
c_{j,k}^{m,1} & \cdots & c_{j,k}^{m,m}
\end{pmatrix}.
\]

\noindent The construction of these matrices relies on the continuity properties of the bilinear \( * \)-product. We can easily check that this vector-valued wavelet respects the VMRA.

\hfill \qed

\section{Construction of wavelet basis for \( L^2(\mathbb{R}^d) \) via tensor products of \( m \)-multiwavelets in \( L^2(\mathbb{R}) \)}\label{meyer}

In \cite{lemarie1986ondelettes}, Lemarié and Meyer constructed a regular orthogonal wavelet basis for \( L^2(\mathbb{R}) \) using properties of the Fourier transform. Extending this one-dimensional framework, they developed multivariate wavelet basis for \( L^2(\mathbb{R}^d) \), \( d \geq 2 \), by employing tensor products of scalar wavelet basis in \( L^2(\mathbb{R}) \), i.e., a scaling function \(\phi\) and a mother wavelet \(\psi\). This pioneering work greatly influenced the development of image processing techniques and became foundational for numerous theoretical studies (see, for example, \cite{jaffard1996wavelet, jaffard2004wavelet, jaffard2007wavelet}). Since then, constructing wavelet basis for \( L^2(\mathbb{R}^d) \) via tensor products has become a standard and effective approach. The novelty of our approach lies in the construction of wavelet basis for \( L^2(\mathbb{R}^d)\) by employing a tensor product of \( m \)-multiwavelets in \( L^2(\mathbb{R}) \) consisting of \( m \)-scaling functions, \( \phi_1, \dots, \phi_m \), and \( m \)-mother wavelets, \( \psi_1, \dots, \psi_m \). This additional flexibility enables the development of multivariate vector-valued wavelets for \( L^2(\mathbb{R}^d, \mathbb{R}^m) \) (see Section \ref{ST}), while preserving key properties from the univariate case, such as orthonormality, compact support, smoothness, and vanishing moments.

\noindent Lemarié and Meyer's construction of a wavelet basis for \( L^2(\mathbb{R}^d) \) was specifically designed for wavelets with compact support in the frequency domain, i.e., wavelets whose Fourier transforms are compactly supported. This construction relied on restrictive sufficient conditions to ensure that the orthonormal family is complete, meaning that it forms a dense subset of \( L^2(\mathbb{R}^d) \) for \( d \geq 2 \). To guarantee that the considered family is complete in \( L^2(\mathbb{R}^d) \), we must take the tensor product of the basis of \( L^2(\mathbb{R}) \), which means forming the tensor product of basis elements overall indices, according to a well-known theorem in functional analysis. In the case of a wavelet basis, in principle, we need to take the product overall scaling indices \( j,k \in \mathbb{Z}^d \). However, by restricting the scaling index to \( j \in \mathbb{Z} \), which is often preferred over \( j \in \mathbb{Z}^d \), their construction seems to select only the diagonal part of the wavelet family. In this sense, to ensure the completeness for \( j \in \mathbb{Z} \), the authors assumed the existence of a subset \( \Omega \subset \mathbb{R}^d \) satisfying the following three conditions:

\begin{equation}
\left\{
\begin{aligned}
    &(1) \quad (\Omega + 2k\pi) \cap (\Omega + 2k'\pi) = \emptyset, \quad \text{for all distinct } k, k' \in \mathbb{Z}^d, \\
    &(2) \quad \text{For every } x \in \mathbb{R}^d, \text{ there exists an integer } j \in \mathbb{N} \text{ such that } 2^j x \in \Omega, \\
    &(3) \quad \text{For each } x \in \mathbb{R}^d, \text{ there exists a vector } \epsilon \in \{ v \in \mathbb{R}^d : v_i = 0 \text{ or } 1,\; i \in \{1, \dots, d\} \} \\
    &\quad \text{such that } \hat{\psi}^{\epsilon_1}(x_1) \neq 0, \dots, \hat{\psi}^{\epsilon_d}(x_d) \neq 0.
\end{aligned}
\right.
\label{total}
\end{equation}

\noindent As Lemarié and Meyer noted, fully satisfying these assumptions is challenging. Their proof was explicitly constructed to approximate these conditions, particularly in the case of \( L^2(\mathbb{R}^2) \) (see \cite{lemarie1986ondelettes} for further details). In reality, their construction is very general and can be applied to all types of wavelets, forming a complete basis without requiring additional assumptions. The key idea is to take the tensor product of the MRA rather than the tensor product of the wavelet basis, as highlighted in \cite{wojtaszczyk1997mathematical}. In this framework, completeness follows directly from the well-known result on the stability of tensor products under increasing unions in Hilbert spaces. In fact, the concept of MRA was introduced by Stéphane Mallat in \cite{mallat1989multiresolution} and appeared in 1989, after their article. In our case, Proposition \ref{HP} follows a procedure similar to previous works based on wavelet-based MRA. However, a key distinction is that we rely on an MRA constructed from \( m \)-multiwavelets rather than traditional wavelets, providing greater flexibility and a richer structure for constructing a wavelet basis for \( L^2(\mathbb{R}^d, \mathbb{R}^m) \), as detailed in Section \ref{ST}.
  
\begin{Proposition}\label{HP}

Let \(\phi = (\phi_1, \dots, \phi_m)\) and \(\psi = (\psi_1, \dots, \psi_m)\) represent an orthonormal \(m\)-multiwavelet basis of \(L^2(\mathbb{R})\). Then, the set of \((2m)^d\) multivariate functions obtained through the tensor product of these components generates an orthonormal wavelet basis for the space \(L^2(\mathbb{R}^d)\). More precisely, the space \( L^2(\mathbb{R}^d) \) admits the following decomposition:
\[
L^2(\mathbb{R}^d) = V_0 \bigoplus_{e=1}^d\bigoplus_{j\in\mathbb{N}} W_j^e,
\]
where
\[
V_0 = \overline{ 
    \operatorname{Vect}\left\{ 
    \phi_{ \alpha_1}(\cdot - k_1) \cdots \phi_{ \alpha_d}(\cdot - k_d)
    : (\alpha_1, \dots, \alpha_d) \in \{1,\dots,m\}^d,\, k \in \mathbb{Z}^d  
    \right\}
}^{\|\cdot\|_{L^2(\mathbb{R}^d)}}
\]

\noindent and  
\[
W_e^j = \overline{\operatorname{Vect} \left\{ \psi_{j, k_1}^{\epsilon_1, \alpha_1}(\cdot) \cdots \psi_{j, k_d}^{\epsilon_d, \alpha_d}(\cdot) : \epsilon \in E_e, \, (\alpha_1, \dots, \alpha_d) \in \{1,\dots,m\}^d, \, k \in \mathbb{Z}^d \right\}}^{\| \cdot \|_{L^2(\mathbb{R}^d)}}.
\]
  
\noindent Here \( E_e \) is the set of vectors in \( \{0, 1\}^d \) with exactly \( e \) nonzero components. The functions \( \psi_{j, k_i}^{0, \alpha_i} \) and \( \psi_{j, k_i}^{1, \alpha_i} \) are defined as:
\[
    \psi_{j, k_i}^{0, \alpha_i}(x) =2^{\frac{j}{2}} \phi_{\alpha_i}(2^{j}x_i - k_i), \quad \psi_{j, k_i}^{1, \alpha_i}(x) =2^{\frac{j}{2}} \psi_{\alpha_i}(2^{j} x_i - k_i).
\]

\end{Proposition}
\begin{Remark}
In particular, according to Corollary \ref{cor}, given any scalar wavelet basis \(\{\phi,\psi\}\) of \( L^2(\mathbb{R} )\), one can construct a $m$-multivariate wavelet basis for \( L^2(\mathbb{R}^d)\).        
\end{Remark}

\noindent \textbf{Proof:} If we consider an \( m \)-scaling function \( \{ \phi_1, \dots, \phi_m \} \) and the corresponding \( m \)-mother wavelets \( \{ \psi_1, \dots, \psi_m \} \) in the space \( L^2(\mathbb{R}) \), according to \cite{chui1996study} we have \( \{ \phi_1, \dots, \phi_m \} \) forms a multiresolution analysis of \( L^2(\mathbb{R}) \) in the form  

$$
\overline{\bigcup_{j \in \mathbb{Z}} V_j} = L^2(\mathbb{R}),
$$

\noindent where the nested subspaces $V_j=\overline{\operatorname{Vect} \{ \phi_i(2^{j}\cdot - k) : i \in \{1, \dots, m\}, k \in \mathbb{Z} \} }$ satisfy  

\[
V_j = V_{j-1} \oplus W_{j-1},
\]

\noindent with the wavelet space given by  

\[
W_{j} = \overline{\operatorname{Vect} \{ \psi_i(2^{j}\cdot - k) : i \in \{1, \dots, m\}, k \in \mathbb{Z} \} }\; \text{for all} \;j\in\mathbb{N}.
\]

\noindent In this proof, we focus on the case \( d=2 \), as the general case follows similarly.
To extend this multiresolution decomposition to \( L^2(\mathbb{R}^2) \), we apply the method of tensor products to the multiresolution decomposition. Specifically, we define the spaces  

\[
\textbf{V}_j = V_j \otimes V_j.
\]

\noindent Since the MRA structure satisfies $
V_j = V_{j-1} \oplus W_{j-1}$, taking tensor products gives  

\[
\textbf{V}_j = (V_{j-1} \oplus W_{j-1}) \otimes (V_{j-1} \oplus W_{j-1}).
\]

\noindent Expanding this expression, we obtain the following.  

$$
\textbf{V}_j = (V_{j-1} \otimes V_{j-1}) \oplus (V_{j-1} \otimes W_{j-1}) \oplus (W_{j-1} \otimes V_{j-1}) \oplus (W_{j-1} \otimes W_{j-1}).
$$

\noindent Using this argument, we deduce that the following sets form orthogonal bases for their respective spaces:

\begin{itemize}
    \item \( \{\phi_{\alpha_1}(2^{j-1}x_1-k_1)\phi_{\alpha_2}(2^{j-1}x_2-k_2) : (\alpha_1,\alpha_2)\in\{1,\dots,m\}^2, (k_1,k_2)\in\mathbb{Z}^2 \} \) is an orthogonal basis for \( V_{j-1} \otimes V_{j-1} \).
    \item \( \{\phi_{\alpha_1}(2^{j-1}x_1-k_1)\psi_{\alpha_2}(2^{j-1}x_2-k_2) : (\alpha_1,\alpha_2)\in\{1,\dots,m\}^2, (k_1,k_2)\in\mathbb{Z}^2 \} \) is an orthogonal basis for \( V_{j-1} \otimes W_{j-1} \).
    \item \( \{\psi_{\alpha_1}(2^{j-1}x_1-k_1)\phi_{\alpha_2}(2^{j-1}x_2-k_2) : (\alpha_1,\alpha_2)\in\{1,\dots,m\}^2,(k_1,k_2)\in\mathbb{Z}^2 \} \) is an orthogonal basis for \( W_{j-1} \otimes V_{j-1} \).
    \item \( \{\psi_{\alpha_1}(2^{j-1}x_1-k_1)\psi_{\alpha_2}(2^{j-1}x_2-k_2) : (\alpha_1,\alpha_2)\in\{1,\dots,m\}^2,(k_1,k_2)\in\mathbb{Z}^2 \} \) is an orthogonal basis for \( W_{j-1} \otimes W_{j-1} \).
\end{itemize}
Following this method, we obtain a multiresolution decomposition of \( L^2(\mathbb{R}^2) \) using \( m \)-multiwavelets:  
\[
L^2(\mathbb{R}^2) = \mathbf{V}_0 \bigoplus_{j\in\mathbb{N}} \left(\mathbf{W}_j^1 \oplus \mathbf{W}_j^2 \oplus \mathbf{W}_j^3 \right),
\]  
where \( \mathbf{W}_j^1 = V_j \otimes W_j \), \( \mathbf{W}_j^2 = W_j \otimes V_j \), and \( \mathbf{W}_j^3 = W_j \otimes W_j \). With a suitable choice of notation, we can consider certain intermediate spaces, such as \( \mathbf{W}_j^1 \oplus \mathbf{W}_j^2 \) as a single space, denoted by \( W_j^{e} \), as mentioned in the proposition. This method can be easily extended to the general case $d\geq 2$. 

\noindent Furthermore, since the \( m \)-scaling functions \( \{\phi_i\}_{i\in\{1,\dots,m\}} \) satisfy the multiresolution analysis conditions, it follows directly that the constructed system of functions  
\[
\{\phi_{ \alpha_1}(\cdot) \cdots \phi_{ \alpha_d}(\cdot)\}_{(\alpha_1, \dots, \alpha_d) \in \{1,\dots,m\}^d}
\]  
forms an MRA of the space \( L^2(\mathbb{R}^d) \). The associated spaces \( \mathbf{V}_j \) are given by  
\[
\mathbf{V}_j = \overline{ 
    \operatorname{Vect}\left\{ 
    \phi_{ \alpha_1}(2^j\cdot - k_1) \cdots \phi_{ \alpha_d}(2^j\cdot - k_d)
    : (\alpha_1, \dots, \alpha_d) \in \{1,\dots,m\}^d,\, k \in \mathbb{Z}^d  
    \right\}
}^{\|\cdot\|_{L^2(\mathbb{R}^d)}}.
\]

\hfill \qed

\section{Construction of multivariate vector-valued wavelet basis for $L^2(\mathbb{R}^d,\mathbb{R}^m)$}\label{ST}

In this section, we present our main contribution, building upon the results of the previous two sections. Specifically, we propose a novel approach for constructing multivariate vector-valued wavelets. Drawing on the insights from \cite{xia1996vector} regarding the relationship between vector-valued and scalar multiresolution analyses, we develop a method to construct an orthonormal vector-valued wavelet basis for \( L^2(\mathbb{R}^d, \mathbb{R}^m) \), derived from wavelet bases of \( L^2(\mathbb{R}, \mathbb{R}^m) \).
 For a given vector-valued wavelet basis \( \phi := (\phi_1, \dots, \phi_m)^T \) and \( \psi := (\psi_1, \dots, \psi_m)^T \) in \( L^2(\mathbb{R}, \mathbb{R}^m) \), we use the connection between \( m \)-multiwavelets and vector-valued multiresolution analysis, as established in Proposition~\ref{XP2}. This allows us to deduce that the functions \( \phi_1, \dots, \phi_m \) and \( \psi_1, \dots, \psi_m \) form an \( m \)-multiwavelet system in \( L^2(\mathbb{R}) \). Using Proposition~\ref{HP}, we then construct an orthonormal basis for \( L^2(\mathbb{R}^d) \). By arranging these scalar-valued functions into a vector, we obtain a vector-valued function. Finally, we get an orthonormal system for \( L^2(\mathbb{R}^d, \mathbb{R}^m) \).
 As an illustration, we provide an example of this construction in the case of \( L^2(\mathbb{R}^2, \mathbb{R}^2) \).

\begin{Theorem}\label{main}
Given any wavelet basis of \( L^2(\mathbb{R}) \) with scaling function \( \phi \) and mother wavelet \( \psi \) respecting the MRA, one can construct a multivariate vector-valued wavelet basis of \( L^2(\mathbb{R}^d, \mathbb{R}^m) \), \(d,m\in\mathbb{N}^{\ast}\), respecting the VMRA. Moreover, the construction preserves key properties of \( \phi \) and \( \psi \), including regularity, vanishing moments, and compact support.  
    
\end{Theorem}
\begin{Remark}
In particular, given any vector-valued wavelet basis of \( L^2(\mathbb{R}, \mathbb{R}^m) \), one can construct a multivariate vector-valued wavelet basis for \( L^2(\mathbb{R}^d, \mathbb{R}^m) \).        
\end{Remark}

\noindent \textbf{Proof:} We will use the same notations as in Proposition \ref{HP}. 

\noindent Let $\phi$ and $\psi$ be a wavelet basis of $L^2(\mathbb{R})$. By Proposition \ref{rm}, first, we construct a vector-valued wavelet in $L^2(\mathbb{R},\mathbb{R}^m)$. Now, using the connection between vector-valued wavelets and $m$-multiwavelets, specifically Proposition \ref{XP2},  we obtain an $m$-multiwavelet system. Finally, using the tensor product construction given in Proposition \ref{HP}, we get a multivariate wavelet system in $L^2(\mathbb{R}^d,\mathbb{R})$. We recall the structure of the subspaces:  

 \[
V_0 = \overline{ 
    \operatorname{Vect}\left\{ 
    \phi_{ \alpha_1}(\cdot - k_1) \cdots \phi_{ \alpha_d}(\cdot - k_d)
    : (\alpha_1, \dots, \alpha_d) \in \{1,\dots,m\}^d,\, k \in \mathbb{Z}^d  
    \right\}
}^{\|\cdot\|_{L^2(\mathbb{R}^d)}}
\]
 \[
    W_j^e = \overline{\operatorname{Vect} \left\{ \psi_{j, k_1}^{\epsilon_1, \alpha_1}(\cdot) \cdots \psi_{j, k_d}^{\epsilon_d, \alpha_d}(\cdot) : \epsilon \in E_e, \, (\alpha_1, \dots, \alpha_d) \in \{1, 2, \dots, m\}^d, \, k \in \mathbb{Z}^d \right\}}^{\| \cdot \|_{L^2(\mathbb{R}^d)}}.
    \]

\noindent The number of elements in the set  
\[
\left\{\psi_{j, k_1}^{\epsilon_1, \alpha_1}(\cdot) \cdots \psi_{j, k_d}^{\epsilon_d, \alpha_d}(\cdot) \right\}_{\epsilon \in E_e,\; \alpha\in \{1, 2, \dots, m\}^d}
\]
is given by  
\[
\frac{d!}{(d-e)!e!}m^d.
\]

\noindent Arranging these $\frac{d!}{(d-e)!e!}m^d$ elements into vectors of size $m$, we obtain $\frac{d!}{(d-e)!e!}m^{d-1}$ vector functions. Let \( f=(f_1,\dots,f_m)^T \) be a function in \( L^2(\mathbb{R}^d,\mathbb{R}^m) \), so that each function $f_i$, $i\in\{1,\dots,m\}$ can be approximated by the system  
\[
\left\{\psi_{j, k_1}^{\epsilon_1, \alpha_1}(\cdot) \cdots \psi_{j, k_d}^{\epsilon_d, \alpha_d}(\cdot) \right\}_{\epsilon \in E_e,\; \alpha\in \{1, 2, \dots, m\}^d}
\]
by Proposition \ref{HP}. Thus, we can write:  
\[
f_i = \sum_{k\in\mathbb{Z}^d}\sum_{\alpha\in \{1, 2, \dots, m\}^d}c_{k}^{ \alpha,i} \psi_{j, k_1}^{0, \alpha_1} \cdots \psi_{j, k_d}^{0, \alpha_d}  
+ \dots  
+ \sum_{k\in\mathbb{Z}^d,\;j\in\mathbb{N}}\sum_{\alpha\in \{1, 2, \dots, m\}^d,\epsilon\in E_d}c_{j,k}^{\alpha,\epsilon, i} \psi_{j, k_1}^{\epsilon_1, \alpha_1} \cdots \psi_{j, k_d}^{\epsilon_d, \alpha_d}.
\]

\noindent Let \( \{A_l\}_{l=1,\dots,m^{d-1}} \) be a partition of \( \{1, 2, \dots, m\}^d \) such that each partition \( A_l \) has a uniform number of elements equal to \( m \). Thus, every \( A_l \) can be expressed as $A_l=\{\alpha_1^l,\dots,\alpha_m^l\},$ where \( \alpha_r^l=(\alpha_{r,1}^l,\dots,\alpha_{r,d}^l) \). 

\noindent Hence, we can rewrite $ f_i $ as:  

$$f_i=\sum_{k\in\mathbb{Z}^d}\sum_{l=1}^{m^{d-1}}\sum_{r=1}^{m}c_{k}^{ \alpha_r^l,i} \psi_{j, k_1}^{0, \alpha_{r,1}^l} \cdots \psi_{j, k_d}^{0, \alpha_{r,d}^l} 
+\dots+\sum_{k\in\mathbb{Z}^d,\;j\in\mathbb{N}}\sum_{\epsilon\in E_d}\sum_{l=1}^{m^{d-1}} \sum_{r=1}^mc_{j,k}^{\alpha_{r}^l, \epsilon, i} \psi_{j, k_1}^{\epsilon_1, \alpha_{r,1}^l} \cdots \psi_{j, k_d}^{\epsilon_d, \alpha_{r,d}^l}.       $$

\noindent Then, the function $f$ can be approximated by:
\[
f=(f_1,\dots,f_m)^T = \sum_{k\in\mathbb{Z}^d}\sum_{l=1}^{m^{d-1}}C_{k}^{ \alpha_l} \Psi_{j, k}^{0, \alpha_l}  
+\sum_{k\in\mathbb{Z}^d,\;j\in\mathbb{N}}\sum_{l=1}^{m^{d-1}} \sum_{\epsilon\in E_1}C_{j,k}^{\alpha_l,\epsilon} \Psi_{j, k}^{\epsilon, \alpha_l}  
+\dots+\sum_{k\in\mathbb{Z}^d,\;j\in\mathbb{N}}\sum_{l=1}^{m^{d-1}} \sum_{\epsilon\in E_d}C_{j,k}^{\alpha_l, \epsilon} \Psi_{j, k}^{\epsilon, \alpha_l}.
\]

The matrix \( C_{j,k}^{\alpha_l,\epsilon, l} \) is given by:  
\[
C_{j,k}^{\alpha_l,\epsilon}=
\begin{pmatrix}
    c_{j,k}^{\alpha_{1}^l,\epsilon, 1}&\dots&c_{j,k}^{\alpha_{m}^l,\epsilon, 1}\\
    \vdots&\ddots&\vdots\\
    c_{j,k}^{\alpha_{1}^l,\epsilon, m}&\dots&c_{j,k}^{\alpha_{m}^l,\epsilon, m}
\end{pmatrix},
\]

and the functions $\Psi_{j, k}^{\epsilon, \alpha_l}$ are:  
\[
\Psi_{j, k}^{\epsilon, \alpha_l} =
\begin{pmatrix}
    \psi_{j, k_1}^{\epsilon_1, \alpha_{1,1}^l} \cdots \psi_{j, k_d}^{\epsilon_d, \alpha_{1,d}^l} \\
    \vdots \\
    \psi_{j, k_1}^{\epsilon_1, \alpha_{m,1}^l} \cdots \psi_{j, k_d}^{\epsilon_d, \alpha_{m,d}^l}
\end{pmatrix}.
\]

\noindent Using the properties of \( \phi \) and \( \psi \), along with Propositions \ref{lem} and \ref{HP}, as well as the use of the partition \( \{A_l\}_{l=1,\dots,m^{d-1}} \) of \( \{1,\dots ,m\}^d \), we can easily show that these vector functions, under dilation and translation, form an orthonormal basis of \( L^2(\mathbb{R}^d,\mathbb{R}^m) \). Concerning the properties of the constructed vector wavelet, they are inherent naturally from the properties of the initial considered wavelet, similar to the standard case of separable wavelets. Furthermore, our method also preserves the conditions of the VMRA.

\hfill \qed

\noindent We illustrate an example in the case $d=m=2$ for clarity, simplifying the construction for the reader.

\begin{Corollary}\label{P}
Let $\phi:= (\phi_1,\phi_2)^T $, $ \psi := (\psi_1,\psi_2)^T$ be an orthonormal vector-valued wavelet basis of $L^2(\mathbb{R},\mathbb{R}^2)$, then $\{\Phi^i_k\}_{i\in\{1,2\},k\in\mathbb{Z}^2}$ , $\{\Psi^i_{j,k}\}_{i\in\{1,\dots ,4\},k\in\mathbb{Z}^2,j\in\mathbb{N} }$, and $\{\Psi^i_{j,k}\}_{i\in\{5,6\},k\in\mathbb{Z}^2,j\in\mathbb{N} }$ are an orthogonal basis of $L^2(\mathbb{R}^2,\mathbb{R}^2)$ where 
\begin{equation*}
     \begin{cases}
        &\Phi^1_k(x,y)=(\phi_1(x-k_1)\phi_1(y-k_2),\phi_2(x-k_1)\phi_2(y-k_2))^T \\
        & \Phi^2_k(x,y)=(\phi_1(x-k_1)\phi_2(y-k_2),\phi_2(x-k_1)\phi_1(y-k_2))^T,
     \end{cases}
\end{equation*}

\begin{equation*}
     \begin{cases}
        &\Psi^1_{j,k}(x,y)=(\phi_1(2^jx-k_1)\psi_1(2^jy-k_2),\phi_2(2^jx-k_1)\psi_2(2^jy-k_2))^T \\
        & \Psi^2_{j,k}(x,y)=(\phi_1(2^jx-k_1)\psi_2(2^jy-k_2),\phi_2(2^jx-k_1)\psi_1(2^jy-k_2))^T\\
        &\Psi^3_{j,k}(x,y)=(\psi_1(2^jx-k_1)\phi_1(2^jy-k_2),\psi_2(2^jx-k_1)\phi_2(2^jy-k_2))^T \\
        & \Psi^4_{j,k}(x,y)=(\psi_1(2^jx-k_1)\phi_2(2^jy-k_2),\psi_2(2^jx-k_1)\phi_1(2^jy-k_2))^T,
\end{cases}
\end{equation*}
and
\begin{equation*}
\begin{cases}
   &\Psi^5_{j,k}(x,y)=(\psi_1(2^{j}x-k_1)\psi_1(2^{j}y-k_2),\psi_2(2^{j}x-k_1)\psi_2(2^{j}y-k_2))^T \\
   & \Psi^6_{j,k}(x,y)=(\psi_1(2^{j}x-k_1)\psi_2(2^{j}y-k_2),\psi_2(2^{j}x-k_1)\psi_1(2^{j}y-k_2))^T.
\end{cases}
\end{equation*}

\end{Corollary}

\noindent\textbf{Proof:} This corollary follows directly from the previous proof of the theorem by taking \( A_1 = \{(1,1), (2,2)\} \) and \( A_2 = \{(1,2), (2,1)\} \). Considering the partitions \( A_1 \) and \( A_2 \) of $\{1,2\}^2$, it applies directly to the formulas cited in the proof.  recall that in this case $E_0=\{(0,0)\}$, $E_1=\{(1,0),(0,1)\}$ and $E_2=\{(1,1)\}$.
so we have 
$$f_i=\sum_{k\in\mathbb{Z}^2}\sum_{l=1}^{2}\sum_{r=1}^{2}c_{k}^{ \alpha_r^l,i} \psi_{j, k_1}^{0, \alpha_{r,1}^l} \psi_{j, k_2}^{0, \alpha_{r,2}^l} 
+\sum_{k\in\mathbb{Z}^2,\;j\in\mathbb{N}}\sum_{\epsilon\in E_1}\sum_{l=1}^{2} \sum_{r=1}^2c_{j,k}^{\alpha_{r}^l, \epsilon, i} \psi_{j, k_1}^{\epsilon_1, \alpha_{r,1}^l} \psi_{j, k_2}^{\epsilon_2, \alpha_{r,2}^l}+\sum_{k\in\mathbb{Z}^d,\;j\in\mathbb{N}}\sum_{l=1}^{2} \sum_{r=1}^2c_{j,k}^{\alpha_{r}^l, (1,1), i} \psi_{j, k_1}^{1, \alpha_{r,1}^l}  \psi_{j, k_2}^{1, \alpha_{r,2}^l}      $$

\noindent and thus,
\[
f=(f_1,f_2)^T = \sum_{k\in\mathbb{Z}^d}\sum_{l=1}^{2}C_{k}^{ \alpha_l} \Psi_{j, k}^{(0,0), \alpha_l}  
+\sum_{k\in\mathbb{Z}^d,\;j\in\mathbb{N}}\sum_{l=1}^{2} \sum_{\epsilon\in E_1}C_{j,k}^{\alpha_l,\epsilon} \Psi_{j, k}^{\epsilon, \alpha_l}  
+\sum_{k\in\mathbb{Z}^d,\;j\in\mathbb{N}}\sum_{l=1}^{2} C_{j,k}^{\alpha_l, (1,1)} \Psi_{j, k}^{(1,1), \alpha_l}.
\]

\noindent The matrix \( C_{j,k}^{\alpha_l,\epsilon} \) is given by:  
\[
C_{j,k}^{\alpha_l,\epsilon}=
\begin{pmatrix}
    c_{j,k}^{\alpha_{1}^l,\epsilon, 1}&c_{j,k}^{\alpha_{2}^l,\epsilon, 1}\\
    c_{j,k}^{\alpha_{1}^l,\epsilon, 2}&c_{j,k}^{\alpha_{2}^l,\epsilon, 2}
\end{pmatrix},
\]

\noindent and the functions $\Psi_{j, k}^{\epsilon, \alpha_l}$ are:  
\[
\Psi_{j, k}^{\epsilon, \alpha_l} =
(\psi_{j, k_1}^{\epsilon_1, \alpha_{1,1}^l} \psi_{j, k_2}^{\epsilon_2, \alpha_{1,2}^l} ,
    \psi_{j, k_1}^{\epsilon_1, \alpha_{2,1}^l} \psi_{j, k_2}^{\epsilon_2, \alpha_{2,2}^l}
)^T
\]

\hfill \qed

\section{Conclusion and Perspectives}
\subsection{Conclusion}
The demand for multivariate vector-valued wavelets has been increasing across various application fields. In particular, their importance has been emphasized in applications such as fluid dynamics in oceanography and aerodynamics, as well as vector-valued signal processing \cite{cui2008m,fowler2002wavelet}. To address this need, we introduce a novel constructive framework for building a wavelet basis for \( L^2(\mathbb{R}^d, \mathbb{R}^m) \) from scalar wavelets while ensuring key mathematical properties such as compact support, regularity, and vanishing moments. In summary, our approach introduces three fundamental innovations: 

\begin{itemize}
    \item We establish a systematic and explicit method for constructing vector-valued wavelet bases in \( L^2(\mathbb{R}, \mathbb{R}^m) \) by directly providing a vector-valued scaling function and a vector-valued mother wavelet from any univariate scalar wavelet basis. These bases inherit key properties while satisfying the VMRA conditions.
    \item Instead of relying on standard wavelet bases of \( L^2(\mathbb{R}) \), we construct multivariate wavelet bases in \( L^2(\mathbb{R}^d, \mathbb{R}) \) with desirable properties using the tensor products of \( m \)-multiwavelets in \( L^2(\mathbb{R}) \), ensuring compatibility with the MRA framework.
    \item We establish the connection between \( m \)-multiwavelets and vector-valued wavelets, enabling a systematic construction of separable multivariate bases in \( L^2(\mathbb{R}^d, \mathbb{R}^m) \) that satisfy VMRA while ensuring that these bases inherit key structural properties.
\end{itemize}
In conclusion, multivariate vector-valued wavelets provide a powerful mathematical framework for signal analysis, particularly in applications where capturing inter-component correlations is essential. Unlike conventional methods that process each component independently, such as multiwavelets, vector-valued wavelets inherently encode dependencies between components, offering a more structured representation of multichannel data. This property is particularly relevant in fields such as medical imaging, fluid dynamics, and remote sensing, where meaningful correlations exist between different signal components. Even deep learning models, particularly convolutional neural networks (CNNs) and transformer-based architectures, have demonstrated significant success in learning complex correlations from multivariate vector-valued signals. However, these models require large amounts of labeled data, are computationally expensive, and often lack interpretability due to their black-box nature. In contrast, vector-valued wavelets could provide an explicit and mathematically grounded approach to capturing inter-component dependencies, potentially serving as a complementary tool for feature extraction and dimensionality reduction in modern signal processing.

\subsection{Perspectives}\label{discussion}

 Several studies have highlighted that even in the case of \( L^2(\mathbb{R}^2,\mathbb{R}
 ) \), wavelet bases constructed using tensor products are not entirely satisfactory from a numerical perspective. For example, in the field of image processing, these bases have been observed to introduce various issues, such as aliasing, defects in image reconstruction, and a lack of geometric adaptability (see \cite{kovacevic1992nonseparable}, \cite{sweldens1998lifting}, \cite{chappelier2005oriented}, \cite{feilner2005orthogonal}, \cite{da2006nonsubsampled}, and \cite{cotronei2019filters}). The pioneering work on constructing multivariate wavelet bases using tensor products is due to Lemarié and Meyer in \cite{lemarie1986ondelettes}. This type of multivariate wavelet basis is called separable, as its basis functions can be written as the product of one-dimensional functions, allowing independent processing along each coordinate axis. In \cite{cohen1993non}, Cohen and Daubechies argued that separable wavelets give special significance to the \( x \) and \( y \) directions and therefore do not cover all orientations, which can be seen as a limitation. Since then, this argument has been frequently cited in the literature (see \cite{kovacevic1992nonseparable}, \cite{sweldens1998lifting}, \cite{chappelier2005oriented}, and \cite{cotronei2019filters}) to explain the shortcomings of separable wavelets and to highlight the benefits of non-separable wavelets, also known as quincunx wavelets. 

Our results open up promising research directions. We aim to refine our construction to overcome the limitations of separable wavelets outlined above. Additionally, in future work, we intend to extend the methods of \cite{cohen1993non} to construct non-separable multivariate vector-valued wavelets in \( L^2(\mathbb{R}^d, \mathbb{R}^m) \). Addressing these challenges could lead to significant advancements in the practical applications of vector-valued wavelets. Furthermore, future research could explore hybrid approaches that integrate vector-valued wavelets with deep learning models to combine the interpretability of wavelet-based methods with the predictive power of deep learning. Such a synergy could provide more robust and interpretable feature extraction techniques for multivariate signal processing, potentially reducing the dependence on large labeled datasets while maintaining strong generalization capabilities. Investigating how vector-valued wavelets can serve as priors for deep learning architectures or as structured representations for unsupervised feature learning could be a promising direction.

\printbibliography

\end{document}